\documentclass[11pt,a4paper]{article}
\usepackage[cp1251]{inputenc}
\usepackage{amsmath}
\usepackage{amssymb}
\usepackage[english]{babel}
\usepackage[dvips]{graphicx}
\usepackage{caption2}
\begin{document}
\title{\textbf{On different curvatures of spheres in Funk geometry}}
\author{Eugene A. Olin}
\maketitle

\begin{center}
\textit{ Geometry Department, Mech.-Math. Faculty, Kharkov National
University, Svoboda sq., 4, 61022-Kharkov, Ukraine.}\\
\textbf{E-mail:} evolin@mail.ru
\end{center}

\begin{abstract}
We compute the series expansions for the normal curvatures of hyperspheres, the Finsler and Rund curvatures of circles in Funk geometry  as the radii tend to infinity. These three curvatures are different at infinity in Funk geometry.
\end{abstract}


\section{Introduction}

A smooth connected manifold $M^n$ is called a \textit{Finsler} manifold (\cite{BCS}) if there is a smooth positively homogeneous on the coordinates in tangent spaces function $F:TM^n\rightarrow [0,\infty)$ such that the symmetric bilinear form
 $\mathbf{g}_y(u,v)=g_{ij}(x,y)u^iv^j:T_x M^n\times T_x M^n\rightarrow
\mathbb{R}$ is positively definite for each pair $(x,y)\in TM^n$, where $g_{ij}(x,y) = \frac{1}{2}[F^2(x,y)]_{y^iy^j}$.

Funk metrics are the non-reversible solutions of Hilbert 4th problem (\cite{Fu}).

Consider a bounded open convex domain $U$ in $\mathbb{R}^n$ with the Euclidean norm $\|\cdot\|$ and let
$ \partial U$ be a $C^{3}$ hypersurface with positive normal
curvatures. For a point
 $x\in U$ and a tangent vector $y \in
T_xU=\mathbb{R}^n$ let $x_0$ be the intersection point of the ray $x+\mathbb{R}_{+}y$ with
$\partial U$. Then the Funk metric is the special Finsler metric defined as following (\cite{Sh}): \begin{equation}\label{eq-funk}\Theta(x,y)=\|y\|\frac{1}{\|x-x_0\|}\end{equation}

Funk metrics are projectively flat non-reversible Finsler metrics of constant negative flag curvature $-\frac{1}{4}$ (\cite{BCS}).
The generalization of Funk metric, the Funk weak metric was constructed in \cite{GT}.

If $U$ is the unit ball then the Funk metric turns into the special projectively flat Randers metric (\cite{Sh}, \cite{Sh1}):
\begin{equation}\label{eq-funkspec}\Theta(x,y)=\frac{\sqrt{|y|^2-(|x|^2|y|^2-<x,y>^2)}+<x,y>}{1-|x|^2}\end{equation}

The Funk metric induces the distance function $d_U(p,q)$ as following.

For given two distinct points $p$ and $q$ in $U$, let $r$ and
 be the intersection point of the half line
$p+\mathbb{R}_{+}(q-p)$ with $\partial
U$. Then the Funk distance (\cite{Fu}):

\begin{equation}\label{eq-funkdist}d_U(p,q)=\ln \frac{||p-r||}{||q-r||}
\end{equation}

It is clear from this formula that Funk metrics are invariant with respect to affine transformations of $\mathbb{R}^n$.

The symmetrization of the Funk metric $ F(x,y)=\frac{1}{2}(\Theta(x,y)+\Theta(x,-y))$ is called the Hilbert metric. Hilbert geometries are the generalizations of the Klein's model of the hyperbolic geometry, Klein's model corresponds to the case when $U$ is the unit ball. Hilbert
geometries are also projectively flat reversible Finsler spaces of constant negative flag curvature $-1$ (\cite{BCS}).

Motivated by the paper by B. Colbois and P. Verovic (\cite{CV}) where it is proved that the unit sphere in the tangent space to the Hilbert geometry approaches the ellipsoid, as the point tend to $ \partial U$, A.A. Borisenko and the author proved in (\cite{Boo5}) that the normal curvatures of hyperspheres, the Rund curvature, and the Finsler curvature of circles in Hilbert geometry tend to 1 as the radii tend to infinity. It is well-known that the normal curvatures of hyperspheres in the hyperbolic space $\mathbb{H}^{n}$ are equal to $\coth(r)$ and tend to 1 as the radius  $r$ tends to infinity.

Let us recall the definitions of curve curvatures in the Finsler geometry. Unlike the Riemannian geometry, in the Finsler geometry there are several definitions of the curvature of a curve.

The normal curvature of a hypersurface in a Finsler space is defined as follows (\cite{Sh}).
Let $\varphi:N\rightarrow M^n$ be a hypersurface in a Finsler manifold $M^n$. A vector $\mathbf{n} \in T_{\varphi(x)}M^n$ is called a normal vector to $N$ at the point  $x \in N$ if
$\mathbf{g}_{\mathbf{n}}(y,\mathbf{n})=0$ for all $y \in T_xN$. Notice that in general non-reversible
case the vector $-n$ is not a normal vector. The \textit{normal curvature} $\mathbf{k}_\mathbf{n}$ at the point $x\in N$ in a direction $y\in T_xN$ is defined as
\begin{equation}\label{eq-curvform}\mathbf{k}_\mathbf{n} = \mathbf{g}_\mathbf{n}(\nabla_{\dot{c}(s)}\dot{c}(s)|_{s=0},\mathbf{n}),\end{equation}
where $\dot{c}(0)=y$, and $c(s)$ is a geodesic in the induced connection on $N$, $\mathbf{n}$ is the chosen unit normal vector.

For a curve $c(s)$ parameterized by its arc length in $M^n$ it is possible to define two more curvatures.

The Finsler curvature of $c(s)$ (\cite{RR}, \cite{F}) is defined as

\begin{equation}\label{eq-curvf}\mathbf{k}_F(c(s)) = \sqrt{\mathbf{g}_{\dot{c}(s)}(\nabla_{\dot{c}(s)}\dot{c}(s),\nabla_{\dot{c}(s)}\dot{c}(s))} \end{equation}

The Rund curvature of $c(s)$ (\cite{RR}) is defined as
\begin{equation}\label{eq-curvr}\mathbf{k}_R(c(s)) = \sqrt{\mathbf{g}_{\nabla_{\dot{c}(s)}\dot{c}(s)}(\nabla_{\dot{c}(s)}\dot{c}(s),\nabla_{\dot{c}(s)}\dot{c}(s))} \end{equation}

Despite the fact that the Funk metric is the metric of constant negative flag curvature it does not necessary possess the same properties as the negatively-curved Riemannian metrics. The non-reversibility of Funk metric leads to the existence of two usually non-collinear normal vector fields and moreover of two families of spheres. So it is a natural question to calculate the different curvatures of spheres in Funk geometry. Since the Funk metric is only positively complete (\cite{Sh}) we will consider only spheres which have arbitrary radii.

Au usual let us call the normal vector which is directed inside (outside) the sphere the inner(outer) normal vector.

\textbf{Theorem 1.} \textit{The normal curvatures of the hyperspheres centred at the same point tend to 0 with respect to the inner normal vector and tend to -2 with respect to the outer normal vector as their radii tend to infinity, uniformly at the point of the hypersphere and in the tangent vector at this point of the hypersphere.}

So the Funk spheres are asymptotically flat with respect to the inner normal vector.

\textbf{Theorem 2.} \textit{The Finsler curvature of the circles centred at the same point tend to 1 as their radii tend to infinity, uniformly at the point of the circle.}

\textbf{Theorem 3.} \textit{The Rund curvature of the circles centred at the same point tend to $\sqrt{2} - 1$ as their radii tend to infinity, uniformly at the point of the circle.}

\section{Spheres in Funk Geometry}

Consider a bounded open convex domain $U \subset \mathbb{R}^{n+1}$
whose boundary is a $C^3$ hypersurface with positive normal
curvatures in $\mathbb{R}^n$.

Fix a point $o \in U$, we will consider this point as the origin
and the center of all the balls. Denote by
$\omega(u):\mathbb{S}^n\rightarrow \mathbb{R}_{+}$ the radial
function for $\partial U$, i. e. the mapping $\omega(u)u$, $u \in
\mathbb{S}^n$ is a parametrization of $\partial U$.

The it is possible to define two families of spheres: forward spheres $S_r^{+} = \{q\in U : d_U(o,q) = r\}$ and backward spheres $S_r^{-} = \{p\in U : d_U(p,0) = r\}$. Denote by $\rho_{r}^{+}(u):\mathbb{S}^n\rightarrow \mathbb{R}_{+}$ and
$\rho_{r}^{-}(u):\mathbb{S}^n\rightarrow \mathbb{R}_{+}$ the radial functions of these two families.

The function $\rho_{r}^{+}(u)$ clearly satisfies
$$r = \ln \frac{\omega(u)}{\omega(u) - \rho_{r}^{+}(u)} $$
therefore we get the explicit formula
\begin{equation}\label{eq-spherefunk} \rho_{r}^{+}(u) = \omega(u) - e^{-r}\omega(u) \end{equation}
From the difference \begin{equation}\label{eq-spherefunkas}\omega(u) - \rho_{r}^{+}(u)= e^{-r}\omega(u)\end{equation} we see that the sphere $S_r^{+}$ tends to the $\partial U$ as $r$ tends to infinity.
But this property does not hold for backward spheres. Indeed, the function $\rho_{r}^{-}(u)$ satisfies
$$r = \ln \frac{\omega(-u) + \rho_{r}^{-}(u)}{\omega(-u)}$$
and hence
$$ \rho_{r}^{-}(u) =  e^{r}\omega(-u) - \omega(-u)$$
Consider the difference
$\omega(u) - \rho_{r}^{-}(u) = \omega(u) -  e^{r}\omega(-u) + \omega(-u) $. The sphere $S_r^{-}$ is contained in $\partial U$ if $\omega(u) - \rho_{r}^{-}(u) > 0$ which is equivalent to
$ r < \ln \left(1 + \frac{\omega(u)}{\omega(-u)}\right)$. Denote by $C_U = \max_{u\in \mathbb{S}^n} \frac{\omega(u)}{\omega(-u)}$ the coefficient of dissymmetry of $U$ with respect to $o$. Then we conclude that the maximal radius of the sphere $S_r^{-}$ which exists in Funk geometry is $\ln \left(1 +C_U \right)$.

In the paper we will consider only forward spheres.

\section{Formulae for Curvatures}

The Chern-Rund covariant derivative along the curve $c(t)$ in the Finsler space equipped with the Funk metric $\Theta$ is given by the formula (\cite{Sh})
\begin{equation}\label{eq-nablacc}\nabla_{c'(t)}c'(t) = \left\{  c''(t)^i+\Theta(c(t),c'(t))c'(t)^i \right\} \frac{\partial}{\partial x^i}\end{equation}

For calculating the normal curvature \eqref{eq-curvform}, Rund curvature \eqref{eq-curvr} and Finsler curvature \eqref{eq-curvf} we need the covariant derivative $\nabla_{\dot{c}(s)}\dot{c}(s)$ of the curve $c(s)$ parameterized by its arc length.

For a given curve $c(t)$ we will denote by the dot the derivative with respect to the arc length $s$, and by the prime the derivative with respect to $t$. Then let $t=t(s)$ be the reparameterization. We get
$$\dot{c}(s)=c'(t)t'_s$$
Using that $s$ in the length parameter we get
$$1=\Theta(c(t),c'(t))t'_s$$
Hence,
$$\dot{c}(s)=\frac{c'(t)}{\Theta(c(t),c'(t))}$$
Next step is to calculate $\nabla_{\dot{c}(s)}\dot{c}(s)$.

$$
\begin{array}{c}
\displaystyle\nabla_{\dot{c}(s)}\dot{c}(s) = \nabla_{\frac{c'(t)}{\Theta(c(t),c'(t))}}\frac{c'(t)}{\Theta(c(t),c'(t))}=\\[10pt]
\displaystyle\frac{1}{\Theta(c(t),c'(t))}\left( \nabla_{c'(t)}\left(\frac{1}{\Theta(c(t),c'(t))} \right)c'(t) + \frac{1}{\Theta(c(t),c'(t))}\nabla_{c'(t)}c'(t)   \right)\\[10pt]
\end{array}
$$
According to (\cite{BCS}),
$$\nabla_{c'(t)}\left(\frac{1}{\Theta(c(t),c'(t))} \right)=-\frac{\mathbf{g}_{c'(t)}(\nabla_{c'(t)}c'(t),c'(t))}{\Theta(c(t),c'(t))^3}$$
Then the derivative $\nabla_{\dot{c}(s)}\dot{c}(s)$ has the form
$$\nabla_{\dot{c}(s)}\dot{c}(s)=\frac{1}{\Theta(c(t),c'(t))^2}\left(\nabla_{c'(t)}c'(t)-\frac{\mathbf{g}_{c'(t)}(\nabla_{c'(t)}c'(t),c'(t))}{\Theta(c(t),c'(t))^2}c'(t)\right)$$
Finally, using \eqref{eq-nablacc} we get the formula:
\begin{equation}\label{eq-nabl}\nabla_{\dot{c}(s)}\dot{c}(s)=\frac{c''(t)+c'(t)\left(\Theta(c(t),c'(t))-\frac{\mathbf{g}_{c'(t)}(\nabla_{c'(t)}c'(t),c'(t))}{\Theta(c(t),c'(t))^2}\right)}{\Theta(c(t),c'(t))^2}\end{equation}

From \eqref{eq-nabl} we get the normal curvature of the curve $c(t)$:
\begin{equation}\label{eq-formnorm}\mathbf{k}_{\mathbf{n}}(t)=\mathbf{g}_{\mathbf{n}}(\nabla_{\dot{c}(t)}\dot{c}(t),\mathbf{n})=\frac{\mathbf{g}_{\mathbf{n}}(c''(t),\mathbf{n})}{\Theta(c(t),c'(t))^2}\end{equation}
Here we used that the normal vector $n$ is by definition $\mathbf{g}_{\mathbf{n}}$-orthogonal to $c'(t)$.

\section{The Choice of the Coordinate System}

Consider the Funk geometry based on a two-dimensional domain $U$ in the Euclidean plane. Fix a point $o$ in the domain $U$ and a point $p \in \partial U$. Since $\partial U$ is a convex curve, it admits the polar representation $\omega(\varphi)$ from the point $o$ such that the point $p$ corresponds to $\varphi = 0$.

Choose the coordinate system on the plane with the origin $O$ at the point $p$; let the axis $x_2$ be orthogonal to $\partial U$ at $p$, $x_1$ be tangent to $\partial U$ at $p$ and $U-\{p\}$ lie in the half-plane $x_2>0$.

In this section we will construct such a affine transformation $P$ of the plane that sends $U$ to $\hat{U}$ and has the following properties:
\begin{enumerate}
\item $P(p) = p$;
\item The vector $u=(0,1)$ is orthogonal to $\partial \hat{U}$ at the point $p$;
\item $\partial \hat{U}$ is the graph of the function $x_2=\hat{f}(x_1)$ such that $\hat{f}(0)=0$, $\hat{f}'(0)=0$, $\hat{f}''(0)=\frac{1}{2}$ in the neighbourhood of  $p$;
    \end{enumerate}

 We are going to give the explicit expression for this transformation and show that after this transformation the curvature of $\partial \hat{U}$ and the derivatives of $\hat{f}$ remain uniformly bounded.

We will use the following lemma that gives the upper bound on the angle between the radial and normal direction to the convex curve.

\textbf{Lemma.} (\cite{Boo2}) \textit{
  Let $\gamma$ be a closed embedded curve in the Euclidean plane whose curvature is greater or equal than $k$. Let $o$ be a point in the interior of the set bounded by $\gamma$, $\omega_0$ be the distance from $o$ to $\gamma$, $\varphi$ be the angle between the outer normal vector at the point $p\in \gamma$ and the vector $op$. Then}

 \begin{equation}
  \label{lem-angle}\cos\angle(u_m,N(m))\geqslant \omega_0 k\end{equation}

Denote by $k$ and $K$ the minimum and maximum of the curvatures of $\partial U$. Denote by $\omega_0=\min\limits_{\varphi}\omega(\varphi)$, $\omega_1=\max\limits_{\varphi}\omega(\varphi)$.

Let the length of the chord of $U$ in the direction $u$ equal $H$, the distance from $o$ to the origin equal $\omega_u$, $\omega_0\leqslant\omega_u\leqslant\omega_1$, and the angle between $u$ and $x_2$ equal $\alpha$.

\textbf{Step 1.} Construct such an affine transformation that makes the vector $\vec{oO}$ parallel to $x_2$. This transformation sends the points $(0,0)$ and $(1,0)$ to themselves, the point $(H \sin \alpha, H \cos \alpha) \in \partial U$ to the point $(0,H)$ and has the expression:
  \begin{equation}
  \label{eq-aff}
           \left\{
           \begin{array}{l}
            \bar{x}_1= x_1-\tan \alpha x_2 \\
             \bar{x}_2=\frac{x_2}{\cos \alpha}\\
                    \end{array}
           \right.
       \end{equation}

Denote the image of $U$ as $\bar{U}$. The point $o$ now has the coordinates $(0,\omega_u)$. Denote by
 $\bar{k}$ the minimum of the curvature of $\partial \bar{U}$ in the $(\bar{x}_1,\bar{x}_2)$ coordinate system, and by $\bar{\omega}_0$ denote the distance from the point $(0,\omega_u)$ to $ \partial \bar{U}$. Note that the eigenvalues of the transformation \eqref{eq-aff} are equal to $1$ and $\frac{1}{\cos \alpha}$, hence
 \begin{equation}
  \label{eq-omega}
    \omega_0 \leqslant \bar{\omega}_0 \leqslant \frac{1}{\cos \alpha}\omega_0
       \end{equation}

Lemma \eqref{lem-angle} than implies that the curvature of $\partial \bar{U}$ remains bounded and separated from zero.

\textbf{Step 2.} Construct such a transformation that the distance from $(0,\omega_u)$ to the origin will equal to 1, and the curvature of $\partial \bar{U}$ at the origin will equal to $1/2$. This transformation has the expression:
 \begin{equation}
  \label{eq-aff2}
           \left\{
           \begin{array}{l}
            \hat{x}_1= \frac{\bar{x}_1}{\omega_u} \\
            \hat{x}_2=\frac{ \bar{x}_2}{2\omega_u^2\bar{k}(0)}\\
                    \end{array}
           \right.
       \end{equation}

Denote the image of $\bar{U}$ as $\hat{U}$. It is obvious that the curvature of  $\partial \hat{U}$ remains bounded.

The announced transformation $P$ is the composition of the transformations \eqref{eq-aff},  \eqref{eq-aff2}, and the following proposition holds:

\textbf{Proposition 1.}\label{th-procurv}\textit{There exists the constant $C_0$ depending on $U$ such that the curvature of $P(\partial \hat{U})$ is bounded from above by $C_0$.}

Let $\partial U$ is the graph of the function $x_2=f(x_1)$ in the initial coordinates system. After the transformation P, $P(\partial U)$ can be considered as the graph of the function $x_2=\hat{f}(x_1)$ such that $\hat{f}(0)=0$, $\hat{f}'(0)=0$, $\hat{f}''(0)=\frac{1}{2}$ in the neighbourhood of $p$.

Finally estimate the third derivative $\hat{f}'''(0)$. Evidently, under the affine transformations \eqref{eq-aff} and \eqref{eq-aff2} the third derivative remains bounded. As $\partial U$ is the compact curve, we obtain

\textbf{Proposition 2.} \textit{There exist the constants $C_1$, $C_2$ depending on $U$, such that $C_1\leqslant\hat{f}'''(0)\leqslant C_2$.}

Analogously we can estimate all higher derivatives.

The Funk metrics for the domains $U$ and $\hat{U}$ are isometric. Therefore without loss of generality we will consider the Funk metric for the domain $\hat{U}$ and will denote $\hat{U}$  by $U$.

\section{Series Expansions for the Metric Tensor of the Funk metric}

Here we will use the method developed in \cite{Boo5}.

Okada lemma (\cite{Ok}) for Funk metrics gives the expression of the derivatives of $\Theta(x,y)$ with respect to the coordinates on tangent spaces through the derivatives with respect to the coordinates on $U$:
$$\Theta(x,y)_{x^k}=\Theta(x,y)\Theta(x,y)_{y^k}$$
Using this lemma we can write:
\begin{equation}
\label{eq-gij}
\begin{array}{c}
\displaystyle g_{ij}(x,y)=\mathrm{\Theta}(x,y)\mathrm{\Theta}_{y^iy^j}(x,y)+\mathrm{\Theta}_{y^i}(x,y)\mathrm{\Theta}_{y^j}(x,y)=\\[10pt]
\displaystyle =\mathrm{\Theta}(x,y)\frac{\Theta_{x^ix^j}(x,y)\Theta(x,y)-2\Theta_{x_i}(x,y)\Theta_{x_j}(x,y)}{\Theta(x,y)^3}+\frac{\Theta_{x^i}(x,y)}{\Theta(x,y)}\frac{\Theta_{x^j}(x,y)}{\Theta(x,y)}\\[10pt]
\end{array}
\end{equation}

For convenience we will use lower indices  $x_i$ for coordinates. Let $\Theta(x_1,x_2,y_1,y_2)$ be the Funk metric. Assume that the point $(x_1,x_2)$ is sufficiently close to $\partial U$. Then we can express $\partial U$ as the graph $x_2=f(x_1)$ such that $f(0)=0$, $f'(0)=0$, $f''(0)=\frac{1}{2}$. Consider a point  $(x_1,x_2)$ above the graph $x_2=f(x_1)$. Denote by $t(x_1,x_2,y_1,y_2)$ the parameter corresponding to the intersection points of the curve $x_2=f(x_1)$ with the line

$$
           \left\{
           \begin{array}{rcl}
            x_1(t) = x_1 + t y_1 \\
             x_2(t) = x_2+ t y_2 \\
           \end{array}
           \right.
$$
 Then

 \begin{equation}\label{eq-theta}\Theta(x_1,x_2,y_1,y_2)=t(x_1,x_2,y_1,y_2)^{-1}\end{equation}

Obtain the derivatives of $\Theta(x_1,x_2,y_1,y_2)$ on $x_1$, $x_2$. The parameter $t(x_1,x_2,y_1,y_2)$
satisfies the functional equation

\begin{equation}\label{eq-main}x_2+t y_2=f(x_1+ t(x_1,x_2,y_1,y_2) y_1)\end{equation}

Differentiate equation \eqref{eq-main} on $x_1$, $x_2$:

\begin{equation}
\label{eq-mg1}
t_{x_1} y_2=f'(x_1+ t y_1)(1+ t_{x_1}y_1), \ 1+t_{x_2} y_2= f'(x_1+ t y_1)t_{x_2}y_1
\end{equation}

We obtain the explicit expressions for $t_{x_1}$, $t_{x_2}$:

\begin{equation}
\label{eq-mg2}
 t_{x_1}=\frac{f'(x_1+ t y_1)}{y_2-y_1f'(x_1+ t y_1)}, \
 t_{x_2}=\frac{1}{y_1f'(x_1+ t y_1)-y_2}
\end{equation}

Differentiating of \eqref{eq-mg1} leads to:

\begin{equation}
\begin{array}{c}
\label{eq-mg11}
\displaystyle y_2 t_{x_1x_1}=f''(x_1+ t y_1)(1+y_1 t_{x_1})^2+f'(x_1+ t y_1)y_1t_{x_1x_1} \\[10pt]
\displaystyle y_2 t_{x_1x_2}=f''(x_1+ t y_1)(1+y_1t_{x_1})y_2t_{x_2}+f'(x_1+ t y_1)y_1t_{x_1x_2} \\[10pt]
\displaystyle y_2 t_{x_2x_2}=f''(x_1+ t y_1)(y_1 t_{x_2})^2+f'(x_1+ t y_1)y_1t_{x_2x_2} \\[10pt]
 \end{array}
\end{equation}

We obtain the expressions for second derivatives of $t$:

\begin{equation}
\begin{array}{c}
\label{eq-mg22}
\displaystyle t_{x_1x_1}=\frac{f''(x_1+ t y_1)(1+y_1t_{x_1})^2}{y_2-y_1f'(x_1+ t y_1)}\\[10pt]
\displaystyle t_{x_1x_2}=\frac{f''(x_1+ t y_1)(1+y_1t_{x_1})y_1t_{x_2}}{y_2-y_1f'(x_1+ t y_1)}\\[10pt]
\displaystyle t_{x_2x_2}=\frac{f''(x_1+ t y_1)(y_1t_{x_2})^2}{y_2-y_1f'(x_1+ t y_1)}\\[10pt]
\end{array}
\end{equation}

Now it is possible to calculate the derivatives of the Funk metric. Formula \eqref{eq-theta} implies
\begin{equation}\label{eq-thetader1}\Theta_{x_k}=-\frac{t_{x_k}}{t^2}\end{equation}

After differentiating \eqref{eq-thetader1} we obtain:
\begin{equation}\label{eq-thetader2}\Theta_{x_kx_l}=-\frac{t_{x_kx_l}t^2-2tt_{x_l}t_{x_k}}{t^4}=(2\Theta^3 t_{x_k}t_{x_l}-\Theta^2 t_{x_kx_l})\end{equation}

Finally, from the formula \eqref{eq-gij} it is possible to obtain the coefficients of the metric tensor.

\subsection{Expansions for $g_{ij}(0,x_2,0,\pm 1)$}

We will need the values of the $g_{ij}(x_1,x_2,0,\pm 1)$ at the points $(x_1,x_2)=(0,x_2)$.

Note that the functions $t(0,x_2,0,\pm 1)$ have the forms
$$ t(0,x_2,0,1)=H-x_2,\,\, t(0,x_2,0,-1)=x_2$$
Here $H$ denotes the length of the chord of $\partial U$ in the direction $(0,1)$.

Consequently
\begin{equation}\label{eq-F21}\mathrm{\Theta}(0,x_2,0,-1)=\frac{1}{x_2}\end{equation}
\begin{equation}\label{eq-F22}\mathrm{\Theta}(0,x_2,0,1)=\frac{1}{H-x_2}\end{equation}
We can estimate the derivatives of the Funk metrics $\Theta(0,x_2,0,\pm1)$.
Formulae \eqref{eq-mg2} imply that at the point $(0,x_2)$

$$t_{x_1}(0,x_2,0,\pm 1)=0,\,\, t_{x_2}(0,x_2,0,\pm 1)=\mp 1$$
$$t_{x_1x_2}(0,x_2,0,\pm 1)=t_{x_2x_2}(0,x_2,0,\pm 1)=0$$

It follows from \eqref{eq-thetader1}, \eqref{eq-thetader2} that
\begin{equation}
\label{eq-thnorm1}
\Theta_{x_2}(0,x_2,0,-1)=\frac{1}{x_2^2}, \ \Theta_{x_2}(0,x_2,0,1)=-\frac{1}{(H-x_2)^2}\\[10pt]
\end{equation}
\begin{equation}
\label{eq-thnorm2}
 \Theta_{x_2x_2}(0,x_2,0,-1)=\frac{2}{x_2^3}, \ \Theta_{x_2x_2}(0,x_2,0,1)=\frac{2}{(H-x_2)^3}\\[10pt]
\end{equation}

Using \eqref{eq-gij}, \eqref{eq-thnorm1}, \eqref{eq-thnorm2}, we get the formulae:
$$g_{12}(0,x_2,0,\pm1)=0$$
\begin{equation}\label{eq-garb}g_{22}(0,x_2,0,-1)=\frac{1}{x_2^2}, g_{22}(0,x_2,0,1)=\left(\frac{1}{H-x_2}\right)^2\end{equation}

\subsection{Expansions for $g_{ij}(0,x_2, 1,0)$}

Here we compute $\Theta$ and $g_{ij}$ at the points $(0,x_2,1,0)$.

 Note that the strict convexity of $\partial U$ implies that $f'(t(x_1,x_2))\neq 0$ for $t(x_1,x_2) \neq 0 $. Then from \eqref{eq-mg2}  we deduce
\begin{equation}\label{eq-t1}t_{x_1}(0,x_2,1,0)=-1\end{equation}
\begin{equation}\label{eq-t2}t_{x_2}(0,x_2, 1,0)=\frac{1}{f'(t(0,x_2, 1,0))}\end{equation}
and from the formulae \eqref{eq-mg22}:
\begin{equation}\label{eq-t11}t_{x_1x_1}(0,x_2, 1,0)=t_{x_1x_2}(0,x_2, 1,0)=0\end{equation}
\begin{equation}\label{eq-t22}t_{x_2x_2}(0,x_2, 1,0)=-\frac{f''(t(0,x_2, 1,0))}{f'(t(0,x_2, 1,0))^3}\end{equation}
Expanding the functional equation \eqref{eq-main} in a power series with respect to $t$ as $x_2\rightarrow0$ we find the expansions of $t(0,x_2,1,0)$
\begin{equation}\label{eq-tayl}x_2=\frac{1}{4}t^2+\frac{1}{6}f'''(0)t^3+O(t^4)\end{equation}
We will find $t$ in expanded form
\begin{equation}\label{eq-t}t=A+B\sqrt{x_2}+Cx_2+Dx_2^{3/2}+O(x_2^2)\end{equation}
After substituting \eqref{eq-t} into \eqref{eq-tayl} and moving all members to the left-hand side we obtain the system of equations
\begin{equation}
\label{eq-tt}
\begin{array}{c}
\displaystyle3A^2+2A^3f'''(0)+(6AB+6A^2f'''(0)B)\sqrt{x_2}+\\[10pt]
\displaystyle+(-12+3B^2+6Af'''(0)B^2+6AC+6A^2f'''(0)C)x_2+\\[10pt]
\displaystyle+(2f'''(0)B^2+6BC+12Af'''(0)BC+6AD+6A^2f'''(0)D)x_2^{3/2}+O(x_2^2)=0\\[10pt]
\end{array}
\end{equation}

Choose coefficients $A$, $B$, $C$, $D$ so that the left side of \eqref{eq-tt} is $O(x_2^2)$.
Equating the coefficients under the powers of $x_2$ to zero we obtain two expansions for $t$ which corresponds to the directions $(1,0)$ and $(-1,0)$.
\begin{equation}\label{eq-tpm}t(0,x_2,\pm 1,0)=\pm 2 \sqrt{x_2}-\frac{4}{3}f'''(0)x_2+O(x_2^2)\end{equation}
In our case we get
 \begin{equation}\label{eq-rpm}t(0,x_2, 1,0)= 2 \sqrt{x_2}- \frac{4}{3}f'''(0)x_2+O(x_2^2)\end{equation}
 Later on all power series will be as $x_2\rightarrow0$. So series expansion for the Funk metric $\Theta$ is:

\begin{equation}\label{eq-F}\Theta(0,x_2,1,0)=\frac{1}{t(0,x_2,1,0)}=\frac{1}{2\sqrt{x_2}}+\frac{f'''(0)}{2}+\frac{2f'''(0)^2}{9}\sqrt{x_2}+O(x_2)\end{equation}

Expand the denominator of \eqref{eq-t2} with respect to $t$:
$$
\begin{array}{c}
\displaystyle t_{x_2}(0,x_2, 1,0)=\frac{1}{f'(t(0,x_2, 1,0)}=\\[10pt]
\displaystyle =\frac{1}{f''(0)t(0,x_2, 1,0)+\frac{1}{2}f'''(0)t(0,x_2, 1,0)^2+O(t(0,x_2, 1,0)^3)}\\[10pt]
\end{array}
$$
Using that $f'(0)=0$, $f''(0)=1/2$, and substituting the value of $t$ from \eqref{eq-tpm} we obtain:
$$t_{x_2}(0,x_2,1,0)=\frac{1}{\frac{1}{2}(2 \sqrt{x_2}-\frac{4}{3}f'''(0)x_2)+\frac{1}{2}f'''(0)(2 \sqrt{x_2}-\frac{4}{3}f'''(0)x_2)^2+ O(x_2^2)}$$

Finally,

\begin{equation}\label{eq-r2}t_{x_2}(0,x_2, 1,0)=\frac{1}{\sqrt{x_2}}-\frac{4f'''(0)}{3}+\frac{40f'''(0)^2}{9}\sqrt{x_2}+O(x_2)\end{equation}

From \eqref{eq-t11} we obtain
\begin{equation}\label{eq-r11}t_{x_1x_1}(0,x_2, 1,0)=t_{x_1x_2}(0,x_2, 1,0)=0\end{equation}
And \eqref{eq-t22} implies
$$t_{x_2x_2}(0,x_2, 1,0)=-\frac{f''(t(0,x_2, 1,0))}{f'(t(0,x_2, 1,0))^3}$$
Expand the numerator and denominator in a series with respect to $t$ and use $f'(0)=0$, $f''(0)=1/2$, \eqref{eq-tpm}:
$$
\begin{array}{c}
\displaystyle t_{x_2x_2}(0,x_2,1,0)=-\frac{\frac{1}{2}+f'''(0)t(0,x_2,1,0)+\frac{1}{2}f^{(4)}(0)t(0,x_2,1,0)^2+O(t^3)}{\left(f''(0)t(0,x_2,1,0)+\frac{1}{2}f'''(0)t(0,x_2,1,0)^2+O(t(0,x_2,1,0)^3)\right)^3}=\\[10pt]
\displaystyle =\frac{-\frac{1}{2}-2f'''(0)\sqrt{x_2}+(\frac{4}{3}f'''(0)^2-4f^{(4)}(0))x_2+\frac{16}{3}f'''(0)f^{(4)}(0)x_2^{3/2}+O(x_2^2)}{x_2^{3/2}+4f'''(0)x_2^2+O(x_2^{5/2})}\\[10pt]
\end{array}
$$

Thus
\begin{equation}\label{eq-r22}t_{x_2x_2}(0,x_2, 1,0)=-\frac{1}{2x_2^{3/2}}-\frac{2f^{(4)}(0)}{ \sqrt{x_2}}+O(1)\end{equation}
From \eqref{eq-thetader1}, \eqref{eq-rpm} we find that
$$\Theta_{x_1}(0,x_2,1,0)=\frac{1}{(2 \sqrt{x_2}- \frac{4}{3}f'''(0)x_2+O(x_2^2))^2}$$
and
\begin{equation}\label{eq-theta1}\Theta_{x_1}(0,x_2,1,0)=\frac{1}{4x_2}+\frac{f'''(0)}{3\sqrt{x_2}}+O(1)\end{equation}

From \eqref{eq-rpm} we deduce
$$\Theta_{x_2}(0,x_2,1,0)= -\frac{\frac{1}{\sqrt{x_2}}-\frac{4f'''(0)}{3}+\frac{40f'''(0)}{9}\sqrt{x_2}+O(x_2)}{(2 \sqrt{x_2}- \frac{4}{3}f'''(0)x_2+O(x_2^2))^2}$$
And we obtain the formula
\begin{equation}\label{eq-theta2}\Theta_{x_2}(0,x_2,1,0)=-\frac{1}{4x_2^{3/2}}-\frac{f'''(0)^2}{\sqrt{x_2}}+O(1)\end{equation}

Using the formulae \eqref{eq-thetader2}, \eqref{eq-rpm}, \eqref{eq-r22} we obtain  the expression for the second derivative of the Funk metric:
$$\Theta_{x_2x_2}(0,x_2,1,0)=\frac{3}{8x_2^{5/2}}+\frac{13f'''(0)^2+3f^{(4)}(0)}{6x_2^{3/2}}+O\left(\frac{1}{x_2}\right)$$

Finally we can estimate the metric coefficients.  From \eqref{eq-theta}, \eqref{eq-thetader2}, \eqref{eq-rpm},  \eqref{eq-theta1}, \eqref{eq-theta2} we get

 \begin{equation}\label{eq-thetabl11}
\begin{array}{c}
\displaystyle \Theta_{x_1x_1}(0,x_2,1,0)\Theta(0,x_2,1,0)-2\Theta_{x_1}(0,x_2,1,0)\Theta_{x_1}(0,x_2,1,0)= 0\\[10pt]
\end{array}
\end{equation}
\begin{equation}\label{eq-thetabl12}\Theta_{x_1x_2}(0,x_2,1,0)\Theta(0,x_2,1,0)-2\Theta_{x_1}(0,x_2,1,0)\Theta_{x_2}(0,x_2,1,0)=0\end{equation}
\begin{equation}
\label{eq-thetabl22}
\begin{array}{c}
\displaystyle\frac{\Theta_{x_2x_2}(0,x_2,1,0)\Theta(0,x_2,1,0)-2\Theta_{x_2}(0,x_2,1,0)\Theta_{x_2}(0,x_2,1,0)}{\Theta(0,x_2,1,0)^3}=\\[10pt]
\displaystyle=\frac{1}{2x_2^{3/2}}+\frac{2f^{(4)}(0)}{\sqrt{x_2}}+O(1)\\[10pt]
\end{array}
\end{equation}
Finally, using \eqref{eq-F}, \eqref{eq-gij}, \eqref{eq-thetabl11}, \eqref{eq-thetabl12}, \eqref{eq-thetabl22} we obtain the series expansions of the metric tensor of Funk metric:

\begin{equation}
\label{eq-g11}
g_{11}(0,x_2, 1,0)=\frac{1}{4x_2}+ \frac{f'''(0)}{3 \sqrt{x_2}} + O(1)
\end{equation}
\begin{equation}
\label{eq-g12}
g_{12}(0,x_2, 1,0)=- \frac{1}{4x_2^{3/2}}- \frac{f'''(0)^2}{ \sqrt{x_2}} + O(1)
\end{equation}

\begin{equation}
\label{eq-g22}
g_{22}(0,x_2,1,0)=\frac{1}{2x_2^{2}}- \frac{f'''(0)}{ 6 x_2^{3/2}} + O(x_2)
\end{equation}

\section{Proof of the Theorems}

As in section 4 fix a point $o$ in the domain $U$ and a point $p \in \partial U$. The curve $\partial U$ admits the polar representation $\omega(\varphi)$ from the point $o$ such that the point $p$ corresponds to $\varphi = 0$; we also assume that $U$ satisfies the conditions 1)-3) from section 4.
Then one can get that $\omega'(0)=0$, $\omega(0)=1$, $\omega''(0)=1/2$.

From \eqref{eq-spherefunk} we get that the circle of radius $r$ admits the parametrization (such that $\varphi = 0$ corresponds to the direction $(0,-1)$).
$$c(\varphi)=\left(( \omega(\varphi) - e^{-r}\omega(\varphi))\sin \varphi,-( \omega(\varphi) - e^{-r}\omega(\varphi))\cos \varphi\right),$$
where $\omega(\varphi)$ is the polar function of $\partial U$.

Then
\begin{equation}\label{eq-c1}c'(0)=(1-e^{-r},0),r\rightarrow\infty\end{equation}
The second derivative:
\begin{equation}\label{eq-c11}c''(0)=\left(0,\frac{1}{2}+\frac{e^{-r}}{2}\right),r\rightarrow\infty\end{equation}
From \eqref{eq-spherefunkas} we get that at the point of the circle the second coordinate
\begin{equation}\label{eq-x2}x_2=e^{-r}\end{equation}

\textbf{Proof of Theorem 1:}
Note that the normal curvature $\mathbf{g}_\mathbf{n}(\nabla_{\dot{c}(s)}\dot{c}(s),\mathbf{n})$ of a hypersurface at the point $x$  depends only on the tangent vector to the curve $c(s)$ at $x$ (\cite{Sh}).
So, in order to obtain the normal curvature of the Funk hypersphere $S_r$ centred at $o$ at the point $p$ in the tangent direction $w$, we consider the normal curvature of the circle $S_r \cap \Pi$ which lies in the plane $\Pi = span (w,\vec{op})$.

From the equality $g_{12}(0,x_2,0,\pm1)=0$ \eqref{eq-garb} it follows that the unit inner normal vector $\mathbf{n_{+}}$ to the circle at $(0,x_2)$ is exactly $\frac{1}{\Theta(0,x_2,0,1)}(0,1)$ and the unit outer normal vector $\mathbf{n_{-}}$ equals to $\frac{1}{\Theta(0,x_2,0,-1)}(0,-1)$.

Finally, taking into account \eqref{eq-formnorm}, \eqref{eq-F21}, \eqref{eq-F22},  \eqref{eq-garb}, \eqref{eq-F}, \eqref{eq-c11}, \eqref{eq-x2}:
\begin{equation}
\label{eq-formnorm1}
\begin{array}{c}
\displaystyle \mathbf{k}_{\mathbf{n_{+}}}(r)=\frac{\frac{1}{2}g_{22}(0,e^{-r},0,1)}{\Theta(0,e^{-r},1-e^{-r},0)^2 \mathrm{F}(0,e^{-r},0,1)}=\\[10pt]
\displaystyle  =\frac{2 e^{-r}}{H} - \frac{8f'''(0) e^{-3r/2}}{3 H} +O(e^{-2r})\\[10pt]
\end{array}
\end{equation}

\begin{equation}
\label{eq-formnorm2}
\begin{array}{c}
\displaystyle \mathbf{k}_{\mathbf{n_{-}}}(r)=-\frac{\frac{1}{2}g_{22}(0,e^{-r},0,-1)}{\Theta(0,e^{-r},1-e^{-r},0)^2 \mathrm{F}(0,e^{-r},0,-1)}=\\[10pt]
\displaystyle  =-2- \frac{8f'''(0) e^{-r/2}}{3 } + \frac{8f'''(0)^2e^{-r}}{9}+O(e^{-2r})\\[10pt]
\end{array}
\end{equation}

If the Euclidean normal curvatures of the hypersurface $\partial U$ are bounded ($k_2\leqslant k_n \leqslant k_1$) then the curvature of the curve  $\partial U' = \partial U \cap \Pi$ is bounded as well. Indeed consider the point $x \in \partial U' \subset \partial U$. Then the curvature $k(x)$ of $\partial U'$ and the normal curvature $k_n(x)$ of $\partial U$ are related as $k(x) = \frac{k_n(x)}{\cos \beta}$. Here $\beta$ in the angle between the radial and normal direction to $\partial U$ at $x$. Using lemma \eqref{lem-angle} we find that $\omega_0 k_2\leqslant \cos \beta \leqslant 1$. Hence the curvature of $\partial U'$ is uniformly bounded for all $y$. Applying propositions 1 and 2 for the Funk geometry based on $U'$ we get the uniformity of the series expansions \eqref{eq-formnorm1}, \eqref{eq-formnorm2} which ends the proof of Theorem 1.

\textbf{Proof of Theorem 2:}
Compute the Finsler curvature of the circles.

Formula \eqref{eq-nablacc} leads to

\begin{equation}\label{eq-nablc1}
\begin{array}{c}
\displaystyle
\nabla_{c'(0)}c'(0) = c''(0)+c'(0)\Theta(c(0),c'(0))=\\[10pt]
\displaystyle =\left(\frac{1}{2\sqrt{e^{-r}}},\frac{1}{2}\right)+O(e^{-r}) \\[10pt]
\end{array}
\end{equation}
Using \eqref{eq-x2}, \eqref{eq-nablc1} we get
$$
\Theta(c(0),c'(0))-\frac{\mathbf{g}_{c'(0)}(\nabla_{c'(0)}c'(0),c'(0))}{\mathrm{\Theta}(c(0),c'(0))^2}= \frac{1}{2 \sqrt{e^{-r}}}-\frac{2 f'''(0)}{3} + O(1) $$

Therefore,
\begin{equation}\label{eq-nablacdot}\nabla_{\dot{c}(0)}\dot{c}(0)=\frac{\left(\frac{1}{2 \sqrt{e^{-r}}}+O(1),\frac{1}{2} + O(e^{-r}) \right)}{\mathrm{\Theta}(c(0),c'(0))^2}\end{equation}

Calculate the Finsler curvature  \eqref{eq-curvf} using the formulae \eqref{eq-x2}, \eqref{eq-nablacdot}.

$$
\begin{array}{c}
\displaystyle
\mathbf{k}_F(r)^2=\mathbf{g}_{\dot{c}(0)}(\nabla_{\dot{c}(0)}\dot{c}(0),\nabla_{\dot{c}(0)}\dot{c}(0))=\\[10pt]
\displaystyle =1- 2f'''(0)e^{-r/2}+O(e^{-r}), \\[10pt]
\end{array}
$$
Proposition 2 provides the uniformity of the obtained expansion, and Theorem 2 follows.

\textbf{Proof of Theorem 3:}

Note that the formula  \eqref{eq-curvr} for the Rund curvature is just the length of the acceleration vector. So for Funk metric

$$ \mathbf{k}_R(r) = \mathrm{\Theta}(c(0), \nabla_{\dot{c}(0)}\dot{c}(0)) $$

Formula  \eqref{eq-nablacdot} leads to

\begin{equation}\label{eq-nablacdot1}\nabla_{\dot{c}(0)}\dot{c}(0) =(2 e^{-r/2} + \frac{16 f'''(0) e^{-r}}{3}, 2 e^{-r} ) + O(e^{-3r/2})\end{equation}

So we need to estimate the value of the Funk metric $\mathrm{\Theta}$ at the point $(e^{-r},0)$ and direction  \eqref{eq-nablacdot1}.

Denote $e^{-r}$ by $x_2$. Equation \eqref{eq-main} leads to

$$x_2 + 2 t x_2 = \frac{1}{4}(2 \sqrt{x_2} t)^2 + \frac{1}{6}f'''(0)(2 \sqrt{x_2} t)^3 + o(2 \sqrt{x_2} t)^3)$$

Expanding the solution of this equation into the series on $x_2$ as $x_2$ tends to 0, we get two solutions

$$t_1(x_2) = 1+\sqrt{2} + \frac{1}{3}(10f'''(0)+7\sqrt{2f'''(0) })\sqrt{x_2}+O(x_2) $$
$$t_2(x_2) = -1+\sqrt{2}+\frac{1}{3}(10f'''(0)-7\sqrt{2f'''(0) })\sqrt{x_2}+O(x_2)  $$

First solution corresponds to the direction $\nabla_{\dot{c}(0)}\dot{c}(0) \dot{c}(0)$ and the second to the direction $-\nabla_{\dot{c}(0)}\dot{c}(0) \dot{c}(0)$.

So from the formula \eqref{eq-theta} we get:

$$\mathrm{\Theta}(c(0), \nabla_{\dot{c}(0)}\dot{c}(0)) = \frac{1}{ t_1(e^{-r})}$$

Expanding this expression leads to:

$$\mathbf{k}_R(r) = \sqrt{2}-1 - \frac{(10 + 7\sqrt{2} )f'''(0)e^{-r/2}}{9+6\sqrt{2}}+O(e^{-r})$$
Proposition 2 provides the uniformity of this expansion which ends the proof of Theorem 3.

 \end{document}